\magnification1200
\input amstex
\documentstyle{amsppt}
\hsize=15.4truecm
\vsize=22.8truecm
\advance\voffset0.6 truecm
\advance\hoffset0.6 truecm
\TagsOnRight
\NoBlackBoxes
\baselineskip0.54cm

\font\be=cmr7

\def\n{\noindent}
\def\vl{\vskip0.4cm}

\def\vb{\vskip0.25cm}
\def\C{\text C\hskip-0.19cm\raise0.03cm\hbox{\be l}\;}
\def\Q{\text Q\hskip-0.19cm\raise0.03cm\hbox{\be l}\;}
\def\R{\text{I\hskip-0.08cm R}}
\def\N{\text{I\hskip-0.08cm N}}

\def\vn{\vskip0.1cm}

\def\v{\,\raise.7ex\hbox{\boxed{ }}\,}

\document
\topmatter
\title{Non-commutative Chern Characters\\  of the  C*-Algebras of spheres\\  and quantum spheres}\endtitle
\author Nguyen Quoc Tho \endauthor
\affil Department of Mathematics, Vinh University, Vietnam \endaffil
\abstract We propose in this paper the construction of non-commutative
Chern characters of the $C^*$-algebras of spheres and quantum spheres. The final
computation gives us a clear relation with the ordinary $\bold Z/(2)$-graded
Chern characters of tori or their normalizers.
\endabstract
\address Department of Mathematics, Vinh University, Vinh City, Vietnam \endaddress
\endtopmatter
\headline={\hss \folio}
\heading Introduction.\endheading
For compact Lie groups the Chern character $ch:\, K^* (G) \otimes \Q
\longrightarrow H_{DR}^* (G; \Q)$ were constructed. In [4]\, - \,[5] we computed the
non-commutative Chern characters of compact Lie group $C^*$-algebras and of
compact quantum groups, which are also homomorphisms from quantum $K$-groups
into entire current periodic cyclic homology of group $C^*$-algebras (resp., of
$C^*$-algebra quantum groups), 
$ch_{C^*}:\, K_* (C^*(G)) \longrightarrow HE_* (C^*(G))$, (resp., $ch_{C^*}:\,
K_* (C^*_\varepsilon (G)) \longrightarrow HE_* (C^*_\varepsilon (G))$). We
obtained also the corresponding algebraic vesion $ch_{\text{alg}}:\, K_* (C^*(G))
\longrightarrow HP_* (C^*(G))$, which coincides with the Fedosov-Cuntz-Quillen
formula for Chern characters [5]. When $A = C_\varepsilon^* (G)$ we first
computed the $K$-groups of $C_\varepsilon^* (G)$ and the $HE_*
(C_\varepsilon^*(G))$. Thereafter we computed the Chern chacracter  
$ch_{C^*}:\, K_* (C^*_\varepsilon (G)) \longrightarrow HE_*
(C^*_\varepsilon (G)$
as an isomorphism modulo torsions.

Using the results from [4] - [5], in this paper we compute the non-commutative
Chern characters $ch_{C^*}:\, K_* (A) \longrightarrow
HE_*(A)$, for two cases $A = C^* (S^n)$, the $C^*$-algebra of spheres and $A =
C_\varepsilon^*(S^n)$, the $C^*$-algebras of quantum spheres. For compact groups
$G = O (n+1)$, the Chern character $ch:\, K_* (S^n) \otimes \Q \longrightarrow
H^*_{DR} (S^n; \Q)$ of the sphere $S^n = O(n+1)/O(n)$ is an isomorphism (see,
[15]). In the paper, we describe two Chern character homomorphisms
\vskip-0.4cm

$$ch_{C^*}:\, K_* (C^*(S^n)) \longrightarrow HE_*(C^*(S^n))$$
and
\vskip-0.6cm

 $$ch_{C^*}:\, K_* (C^*_\varepsilon (S^n)) \longrightarrow
HE_*(C^*_\varepsilon (S^n)).$$
Finally, we show that there is a commutative diagram 
$$\matrix K_* \big(C^* (S^n)\big) & \overset{ch_{C^*}} \to
{\text{-----}\hskip-0.2cm 
\longrightarrow } &
HE_* \big(C^*(S^n)\big)\\
\bigg\downarrow \cong && \bigg\downarrow \cong \\
K_* \big(\C(\Cal N_{\bold T_n})\big) & \overset {ch_{CQ}}
\to{\text{-----}\hskip-0.2cm \longrightarrow } & 
HE_* \big(\C(\Cal N_{\bold T_n})\big)\\ 
\bigg\downarrow \cong && \bigg\downarrow \cong \\
K^* \big(\Cal N_{\bold T_n}\big) & \overset{ch}\to {\text{-----}\hskip-0.2cm
\longrightarrow} & H_{DR}^* 
\big(\Cal N_{\bold T_n}\big)\endmatrix $$ 
(Similarly, for $A = C_\varepsilon^*(S^n)$, we have an analogous commutative
diagram with $W \times S^1$ of place of $W \times S^n$), from which we 
deduce that $ch_{C^*}$ is an isomorphism modulo torsions. 

We now briefly review the structure of the paper. In section 1, we compute the
Chern chatacter of the $C^*$-algebras of spheres. The computation of Chern
character of $C^*(S^n)$ is based in two crucial points: 

i) Because the sphere $S^n = O(n+1)\big/O(n)$ is a homogeneous space and
$C^*$-algebra of $S^n$ is the transformation group $C^*$-algebra, following
J. Parker [10], we have, $C^*(S^n) \cong C^* (O(n)) \otimes \Cal K (L^2 (S^n)).$ 

ii) Using the stability property theorem $K_*$ and $HE_*$ in [5], we reduce it
to the computation of $C^*$-algebras of subgroup $O(n)$ in $O(n+1)$ group. 

In section 2, we compute the Chern character of $C^*$-algebras of quantum
spheres. For quantum sphere $S^n$, we define the compact quantum
$C^*$-algebra $C_\varepsilon^* (S^n)$, where $\varepsilon$ is a positive real number.
Thereafter, we prove that 
$$C_\varepsilon^* (S^n) \cong \C(S^1) \oplus
\bigoplus\limits_{e \ne \omega \in W} \int_{S^1}^{\oplus} \Cal K (H_{w, t})
dt,$$ 
where $\Cal K (H_{w, t})$ is the elementary algebra of compact operators in a
separable infinite dimensional Hilbert space $H_{w, t}$ and $W$ is the Weyl of a
maximal torus $\bold T_n$ in $SO(n)$. 

Similar to Section 1, we first compute the $K_*(C_\varepsilon^*(S^n))$ and
$HE_*(C_\varepsilon^*(S^n))$, and we prove that $ch_{C^*}:\, K_*
(C_\varepsilon^*(S^n)) \longrightarrow HE_* (C_\varepsilon^*(S^n))$ is a
isomorphism modulo torsions.

\n {\bf Notes on Notation:} For any compact space $X$, we write $K^* (X)$ for the
$\bold Z/(2)$-graded topological $K$-theory of $X$. We use Swan's theorem to
identify $K^*(X)$ with $\bold Z/(2)$-graded $K_*(\C(X))$. For any involutive Banach
algebra $A, K_*(A)$, $HE_*(A)$, $HP_* (A)$ are $\bold Z/(2)$- graded algebraic or
topological $K$-groups of $A$, entire cyclic homology, and periodic cyclic
homology of $A$, respectively. If $\bold T$ is a maximal torus of a compact group
$G$, with the corresponding Weyl group $W$, write $\C(\bold T)$ for the algebra of
complex valued functions on $\bold T$. We use the standard notations from the root
theory such as $P$, $P^+$ for the positive highest weights, etc... We denote by
$\Cal N_{\bold T}$ the normalizer of $\bold T$ in $G$, by $\N$ the set of
natural numbers, $\R$ the fied of real numbers and $\C$ the field of complex
numbers, $\ell_A^2 (\N)$ the standard $\ell^2$ space of square integrable
sequences of elements from $A$, and finally by $C_\varepsilon^*(G)$ we denote
the compact quantum algebras, $C^*(G)$ the $C^*$-algebra of $G$. 
\vb

\n \S {\bf 1. Non-commutative Chern characters of $\bold C^*$-algebras of
spheres.} 
\vb

In this section, we compute non-commutative Chern characters of $C^*$-algebras
of spheres. Let $A$ be an involutive Banach algebra. We construct the
non-commutative Chern characters $ch_{C^*}:\, K_* (A) \longrightarrow HE_* (A)$,
and show in [4] that for $C^*$-algebra $C^*(G)$ of compact Lie groups $G$, the
Chern character $ch_{C^*}$ is an isomorphism. 
\vb

\n {\bf Proposition 1.1} ([5], Theorem 2.6). {\it Let $H$ be a separable
Hilbert space and $B$ an arbitrary Banach space. We have

i)
\vskip-1.5cm

 $$\align & K_*(\Cal K (H)) \cong K_*(\text{\C})\\
& K_* (B\otimes \Cal K (H)) \cong K_*(B); \endalign $$ 

ii) 
\vskip-1.6cm

$$\align & HE_*(\Cal K (H)) \cong HE_*(\text{\C})\\
& HE_* (B\otimes \Cal K (H)) \cong HE_*(B), \endalign$$ 
where $\Cal K(H)$ is the elementary algebra of compact operators in a separable
infitite-dimensional Hilbert space $H$. } \hfill \v

\n {\bf Proposition 1.2.} ([5], Theorem 3.1). {\it Let $A$ be an involutive
Banach algebra with unity. There is a Chern character homomorphism
\vskip-0.3cm

$$ch_{C^*}:\, K_*(A) \longrightarrow HE_* (A).$$} 
\vskip-0.5cm\hfill\v

\n {\bf Proposition 1.3.} ([5], Theorem 3.2). {\it Let $G$ be an compact group
and $\bold T$ a fixed maximal torus of $G$ with Weyl group $W:= \Cal N_{\bold
T}/\bold T$. Then
the Chern character $ch_{C^*}:\, K_* (C^*(G)) \longrightarrow HE_* (C^*(G))$ is
an isomorphism modulo torsions, i.e. 
$$ch_{C^*}:\, K_*(C^*(G)) \otimes \text{\C} \overset \cong \to \longrightarrow HE_*
(C^*(G)),$$ 
which can be identified with the classical Chern character
$$ch:\, K_* (C(\Cal N_{\bold T})) \longrightarrow HE_*
(C(\Cal N_{\bold T})),$$ 
that is also an isomorphism modulo torsions, i.e. }
$$ch:\, K^*(\Cal N_{\bold T}) \otimes \C \overset \cong \to \longrightarrow H_{DR}^*
(\Cal N_{\bold T}).$$\vskip-0.6cm \hfill \v

Now, for $S^n = O(n+1)/O(n)$, where $O(n)$, $O(n+1)$ are the orthogonal
matrix groups. We denote by $\bold T_n$ a fixed maximal torus of $O(n)$ and
$\Cal N_{\bold T_n}$ the normalizer of $\bold T_n$ in $O(n)$. Following
Proposition 1.2, there a natural Chern character $ch_{C^*}:\, K_* (C^*(S^n))
\longrightarrow HE_*(C^*(S^n))$. Now, we compute first $K_*(C^*(S^n))$ and then
$HE_*(C^*(S^n))$ of $C^*$-algebra of the sphere $S^n$. 
\vb

\n {\bf Proposition 1.4.} 
$$HE_*(C^* (S^n)) \cong H_{DR}^W (\bold T_n)$$ 

\n {\it Proof}: We have 
\vskip-0.9cm

$$\align & HE_* (C^*(S^n)) = HE_* \big(C^*(O(n+1)\big/O(n))\big)\\
& \cong HE_* (C^* (O(n))\otimes \Cal K
\big(L^2\big(O(n+1)/O(n)))\big)\endalign$$ 
(in virtue of [10], the $\Cal K\big(L^2\big( O(n+1)\big/O(n))\big)$ is a
$C^*$-algebra compact operators in a separable Hilbert space
$L^2\big(O(n+1)\big/O(n))\big)$
$$\align & \cong HE_* (C^*(O(n)))\qquad \quad\;\; (\text{by Proposition 1.1})\\
&\cong HE_* (\C(\Cal N_{\bold T_n}))\qquad\qquad \quad\quad\;\; (\text{see
[5]}).\endalign $$
Thus, we have $HE_* (C^*(S^n)) \cong HE_* (\C(\Cal N_{\bold T_n}))$. 

Apart from that, because $\C(\Cal N_{\bold T_n})$ is the commutative
$C^*$-algebra, by a result Cuntz-Quillen's [1], we have an isomorphism
$$HP_*(\C(\Cal N_{\bold T_n})) \cong H_{DR}^*(\Cal N_{\bold T_n}).$$
Moreover, by a result of Khalkhali [8] - [9], we have 
$$HP_* (\C(\Cal N_{\bold T_n})) \cong HE_*(\C(\Cal N_{\bold T_n})).$$ 
We have, hence
$$\align HE_*(C^*(S^n)) & \cong HE_* (\C(\Cal N_{\bold T_n})) \cong HP_* (\C(\Cal
N_{\bold T_n}))\\
& \cong H_{DR}^* (\Cal N_{\bold T_n}) \cong H_{DR}^W (\bold T_n)\;\qquad \qquad
(\text{by}\;\, [15]).\endalign$$
\vskip-0.7cm \hfill \v
\vb

\n {\it Remark 1.} Because $H_{DR}^W (\bold T_n)$ is the de Rham cohomology
of 
$\bold T_n$, invariant under the action of the Weyl group $W$, following Watanabe
[15], we have a canonical isomorphism 
$ H_{DR}^W (\bold T_n)  \cong H^* (SO(n))
 = \Lambda_{\text \C} (x_3, x_7,..., x_{2i+3}),$
where $x_{2i+3} = \sigma^*(p_i) \in H^{2n+3} (SO(n))$ and $\sigma^*: H^*
(BSO(n), R)$, 
$\longrightarrow H^* (SO(n), R)$ for a commutative ring $R$ with a unit $1 \in
R$, and $p_i = \sigma_i (t_1^2, t_2^2,..., t_i^2) \in H^*(B\bold T_n, \bold Z)$
the Pontryagin classes. 

Thus, we have 
$$HE_* (C^* (S^n)) \cong \Lambda_{\text \C} (x_3, x_7,..., x_{2i + 3}).$$ 

\n {\bf Proposition 1.5.} 
$$K_* (C^* (S^n)) \cong K^* (\Cal N_{\bold T_n}).$$ 

\n {\it Proof.} We have \vskip-0.7cm 
$$\align & K_* \big(C^* (S^n)\big) = K_* \big(C^* (O(n+1)\big/O(n))\big)\\
& \cong K_* \big(C^*(O(n)) \otimes \Cal K
\big(L^2(O(n+1)\big/O(n))\big)\big) \;\; (\text{see} \;[10])\\
& \cong K_* \big(C^*(O(n))) \qquad\qquad\;\;\; (\text{by Proposition 1.1})\\
& \cong K_* \big(\C(\Cal N_{\bold T_n}))\\
& \cong K^* \big(\Cal N_{\bold T_n})\qquad\qquad\;\;\ (\text{by Lemma 3.3, from}\;
[5]). \endalign$$ 
Thus,\qquad\quad  $K_* (C^*(S^n)) \cong K^* (\Cal N_{\bold T_n}).$ \hfill \v
\vb

\n {\it Remark 2.} Following Lemma 4.2 from [5], we have
$$\align K^* (\Cal N_{\bold T_n}) & \cong K^* \big(SO(n+1)\big)\big/\text{Torsion}\\
& = \Lambda_{\bold Z} (\beta (\lambda_1),..., \beta (\lambda_{n-3}),
\varepsilon_{n+1}\big),\endalign $$ 
where $\beta:\, R (SO(n)) \longrightarrow \widetilde K^{-1} (SO(n))$ be the
homomorphism of Abelian groups assigning to each representation $\rho:\,
SO(n) \longrightarrow U (n+1)$ the homotopy class $\beta(\rho) = [i_n \rho]
\in [SO(n), U] = \widetilde K^{-1} (SO(n))$, where $i_n:\, U(n+1) \to U$ is the
canonical one, $U(n+1)$ and $U$ be the $n\,-\,th$ and infinite unitary groups
respectively and $\varepsilon_{n+1} \in K^{-1}(SO(n+1))$. We have, finally
$$K_* (C^* (S^n)) \cong \Lambda_{\bold Z} (\beta (\lambda_1),...,
\beta(\lambda_{n-3}), \varepsilon_{n+1}).$$ 

Moreover, the Chern character of $S U (n+1)$ was computed in [14], for all $n
\ge 1$. Let us recall the result. Define a function
$$\phi:\, \N \times \N \times \N \longrightarrow \bold Z,$$ given by 
\vskip-0.5cm

$$\phi(n, k, q) = \sum_{i=1}^k (-1)^{i-1} \pmatrix n\\
k-1\endpmatrix i^{q-1}.$$ 
\vb

\n {\bf Theorem 1.6.} {\it Let $\bold T_n$ be a fixed maximal torus of $O(n)$
and $\bold T$ the fixed maximal torus of $SO(n)$, with Weyl groups $W:=
\Cal N_{\bold T}/\bold T$, the Chern character of $C^*(S^n)$ 

$$ch_{C^*} :\,K_*(C^* (S^n)) \longrightarrow HE_* (C^* (S^n))$$ 
is a isomorphism, given by } 
$$\align ch_{C^*} \big(\beta (\lambda_k)\big) & = \sum_{i=1}^n ((-1)^{i-1}
2\big/(2i-1)!) 
\phi (2n+1, k, 2i) x_{2i+3}, \;\, (k = 1, 2,..., n-1)\\
ch_{C^*} (\varepsilon_{n+1}) & = \sum_{i=1}^n \big((-1)^{i-1} 2\big/(2i-1)!)
((1/2^n) \sum_{k=1}^n \phi (2n+1, k, 2i)) x_{2i+3}. \endalign$$
\vb

\n {\it Proof.} By Proposition 1.5, we have 

$$K_* (C^*(S^n)) \cong K_* (\C (\Cal N_{\bold T_n})) \cong K^* (\Cal N_{\bold T_n})$$
and 
\vskip-0.5cm
$$ HE_* (C^*(S^n)) \cong HE_* (\C(\Cal N_{\bold T_n}))
 \cong H_{DR}^* (\Cal N_{\bold T_n}) \qquad  (\text{by Proposition
1.4}).$$

\n Now, consider the commutative diagram 
$$\matrix K_* \big(C^* (S^n)\big) & \overset{ch_{C^*}} \to
{\text{-----}\hskip-0.2cm 
\longrightarrow } &
HE_* \big(C^*(S^n)\big)\\
\bigg\downarrow \cong && \bigg\downarrow \cong \\
K_* \big(\C(\Cal N_{\bold T_n})\big) & \overset {ch_{CQ}}
\to{\text{-----}\hskip-0.2cm \longrightarrow } & 
HE_* \big(\C(\Cal N_{\bold T_n})\big)\\ 
\bigg\downarrow \cong && \bigg\downarrow \cong \\
K^* \big(\Cal N_{\bold T_n}\big) & \overset{ch}\to {\text{-----}\hskip-0.2cm
\longrightarrow} & H_{DR}^* 
\big(\Cal N_{\bold T_n}\big)\endmatrix $$ 
Moreover, by the results of Watanabe [15], the Chern character $ch:\,K^*
\big(\Cal N_{\bold
T_n}\big) \otimes \C \longrightarrow H_{DR}^* \big(\Cal N_{\bold T_n}\big)$ is an
isomorphism. 

Thus,\; $ch_{C^*}:\, K_*\big(C^* (S^n)\big) \longrightarrow HE_*
\big(C^*(S^n)\big)$ is an isomorphism (by Proposition 1.4 and 1.5), given by 
\vskip-0.5cm

$$\align ch_{C^*}(\beta(\lambda_k))  & = \sum_{i=1}^n \big((-1)^{i-1}\;
2\big/(2i-1)!\big) \phi (2n+1, k, 2i) x_{2i+3},\;\, (k = 1, 2,..., n-1),\\
ch_{C^*}(\varepsilon_{n+1}) & = \sum_{i=1}^n \big((-1)^{i-1}
2\big/(2i-1)!\big) \Big((1/2^n) \sum_{k=1}^n \phi (2n+1, k, 2i)\Big)
x_{2i+3},\endalign$$ 
where:
\vskip-0.8cm

$$\align K_* \big(C^* (S^n)\big) & \cong \Lambda_{\bold Z}
(\beta(\lambda_1),..., \beta(\lambda_{n-3}), \varepsilon_{n+1}),\\
HE_*\big(C^* (S^n)\big) & \cong \Lambda_{\C}
(x_3, x_7,..., x_{2n+3}).\endalign$$ \vskip-0.7cm \hfill\v
\vb

\n \S {\bf 2. Non-Commutative Chern character of $\bold C^*$-algebra of quantum
spheres.} 
\vb

In this section, we at first recall definitions and main properties of compact
quantum spheres and their representations. More precisely, for $S^n$, we define
$C_\varepsilon^* (S^n)$, the $C^*$-algebras of compact quantum spheres as the
$C^*$-completion of the *-algebra $\Cal F_\varepsilon (S^n)$ with respect to
the $C^*$-norm, where $\Cal F_\varepsilon(S^n)$ is the quantized Hopf
subalgebra of the Hopf algebra, dual to the quantized universal enveloping
algebra $U(\Cal G)$, generated by matrix elements of the $U(\Cal G)$ modules of type {\bf
1}
(see [3]). We prove that 
$$C_\varepsilon^* (S^n) \cong \C(S^1) \oplus
\bigoplus\limits_{e \ne \omega\in W} \int_{S^1}^{\oplus} \Cal K (H_{w, t})
dt,$$ 
where $\Cal K (H_{w, t})$ is the elementary algebra of compact operators in
a separable infinite-dimensional Hilbert space $H_{w, t}$ and $W$ is the Weyl
group of $S^n$ with respect to a maxinal torus $\bold T$. 

After that, we first compute the $K$-groups $K_*\big(C_\varepsilon^*(S^n)\big)$
and the $HE_* (C_\varepsilon^* 
(S^n))$, respectively. Thereafter we define the Chern character of
$C^*$-algebras quantum spheres, as a homomorphism from $K_* \big(C_\varepsilon^*
(S^n)\big)$ to 
$HE_* \big(C_\varepsilon^* (S^n)\big)$, and we prove that $ch_{C^*}:\, K_*
(C_\varepsilon^* (S^n)) 
\longrightarrow HE_*\big(C_\varepsilon^* (S^n)\big)$ is an isomorphism modulo
torsions.  

Let $G$ be a complex algebraic group with Lie algebra $\Cal G$= Lie $G$ and
$\varepsilon$ is real number, $\varepsilon \ne -1$. 
\vb

\n {\bf Definition 2.1.} ([3], Definition 13.1). {\it The quantized function
algebra $\Cal F_\varepsilon (G)$ is the subalgebra of the Hopf algebra dual to
$U_\varepsilon (\Cal G)$, generated by the matrix elements of the
finte-dimensional $U_\varepsilon (\Cal G)$-modules of type {\bf 1}. }

For compact quantum groups the unitary representation of $\Cal F_\varepsilon
(G)$ are parameterized by pairs $(w, t)$, where $t$ is an element of a fixed
maximal torus of the compact real form of $G$ and $w$ is an element of the Weyl
group $W$ of $\bold T$ in $G$. 

Let $\lambda \in P^+, V_\varepsilon (\lambda)$ be the irreducible $U_\varepsilon
(\Cal G)$-module of type {\bf 1} with the highest weight $\lambda$. Then
$V_\varepsilon (\lambda)$ admits a positive definite hermitian form (.,.), such
that $(xv_1, v_2) = (v_1, x^* v_2)$ for all $v_1, v_2 \in V_\varepsilon
(\lambda), x \in U(\Cal G)$. Let $\{v^\nu_\mu\}$ be an orthogonal basis for weight
space $V_\varepsilon (\lambda)_\mu$, $\mu \in P^+$. Then $\cup\{v_\mu^\nu\}$
is an orthogonal basis for $V_\varepsilon (\lambda)$. Let $C_{\nu, s; \mu,
r}^\lambda (x) = \big(x v_\mu^r, v_\nu^s\big)$ be the associated matrix
elements of $V_\varepsilon (\lambda)$. Then the matrix elements $C_{\nu, s; \mu,
r}^\lambda$ (where $\lambda$ runs
throngh $P^+$, while $(\mu, r)$ and $(\nu, s)$ runs independently through the
index set of a basis of $V_\varepsilon (\lambda)$ form a basis of $\Cal
F_\varepsilon (G)$ (see [3]). 

Now very irreducible * - representation of $\Cal F_\varepsilon (SL_2 (\C))$ is
equivalen to a representation belonging to one of the following two families,
each of which is parameterized by $S^1 = \{ t \in \C \vert \; \vert t\vert = 1\}$, 

i) the family of one-dimensional representation $\Cal T_t$

ii) the family $\pi_t$ of representation in $\ell^2 (\N)$ (see [3])

Moreover, there exists a surjective homomorphism $\Cal F_\varepsilon (G)
\longrightarrow \Cal F_\varepsilon (SL_2(\C))$ induced by the natural inclusion
$SL_2 (\C) \hookrightarrow G$ and by composing the representation $\pi_{-1}$
of $\Cal F(SL_2 (\C))$ with this homomophis, we obtain a representation of $\Cal
F_\varepsilon (G)$ in $\ell^2 (\N)$ denoted by $\pi_{s_i}$, where $s_i$ appears
in the reduced decomposition $w = s_{i_1}. s_{i_2}\ldots s_{i_k}$. More
precisely, $\pi_{s_i}:\, \Cal F_\varepsilon (G) \longrightarrow \Cal L (\ell^2
(\N))$ is of class CCR (see [11]), i.e its image is dense in the ideal of
compact operators in $\Cal L (\ell^2 (\N))$. 

The representation $\Cal T_t$ is one-dimensional and is of the form $$\Cal T_t
\big(C_{\nu, s, \mu, r}^\lambda\big) = \delta_{r, s} \delta_{\mu \nu} exp\,
(2\pi \sqrt{-1} \mu (x)),$$ if $t=\, exp\, (2\pi \sqrt{-1} x) \in \bold T$,
for $x \in$ Lie $\bold T$, (see [3]). 
\vb

\n {\bf Proposition 2.2} ([3], 13.1.7). {\it Every irreducible unitary
representation of $\Cal F_\varepsilon (G)$ on a separable Hilbert space is the
completion of a unitarizable highest weight representation. Moreover, two such
representation are equivalent if and only if they have the same highest
weight.} \hfill\v
\vb

\n {\bf Proposition 2.3.} ([3], 13.1.9). {\it Let $\omega =
s_{i_1}.s_{i_2}\ldots s_{i_k}$ be a resuced decomposition of an element $w$ of
the Weyl group $W$ of $G$. Then 

i) the Hilbert space tensor product $\rho_{w, t} = \pi_{s_{i_1}} \otimes
\pi_{s_{i_2}} \otimes \cdots \otimes \pi_{s_{i_k}} \otimes \Cal T_t$ is an
irreducible 
*-representation of $\Cal F_{\varepsilon} (G)$ which is associated to the
Schubert cell $S_w$; 

ii) up to equivalence, the representation $\rho_{w, t}$ does not depend on the
choice of the reduced decomposition of $w$; 

iii) every irreducible *-representation of $\Cal F_\varepsilon (G)$ is
equivalent to some $\rho_{w, t}$.} \hfill \v

The sphere $S^n$, can be realized as the orbit under the action of the compact
group $S U (n+1)$ of the highest weight vector $v_0$ in its natural
$(n+1)$-dimensional representation $V^{\hbar}$ of $S U (n+1)$. If $t_{rs}, 0 \le
r, s, \le n$, are the matrix entries of $V^{\hbar}$, the algebra of functions on
the orbit is generated by the entries in the ``first column" $t_{s0}$ and their
complex conjugates. In fact, 
$$\Cal F(S^n):= \C [t_{00},\ldots, t_{n0}, \overline t_{00}, \ldots \overline
t_{n0}]\big/\sim,$$  
where ``$\sim$" is the following equivalence ralation 
$$t_{s0} \overline t_{s0} \Longleftrightarrow \sum\limits_{s=0}^n t_{s0}.
\overline t_{s0} = 1.$$
\vn

\n {\bf Proposition 2.4}. ([3], 13.2.6). {\it The *-structure on Hopf algebra
$\Cal F_\varepsilon (SL_{n+1} (\C))$, is given by 
$$t_{rs}^* = (-\varepsilon)^{r-s}. q\det\; (\widehat T_{rs}),$$ 
where $\widehat T_{rs}$ is the matrix obtained by removing the $r^{th}$ row and
the $s^{th}$ 
column from $T$. }
\vb

\n {\bf Definition 2.5}. ([3], 13.2.7). {\it The *-subalgebra of $\Cal
F_\varepsilon (SL_{n+1} (\C))$ generated by the elements $t_{s0}$ and
$t_{s0}^*$, for $s = 0,\ldots n$, is called the quantized algebra of functions
on the sphere $S^n$, and is denoted by $\Cal F_\varepsilon (S^n)$. It is a
quantum $SL_{n+1}$ $(\text{\C})$-space. }

We set $z_s = t_{s0}$ from now on. Using Proposition 2.4, is is easy to see
that the following relations hold in $\Cal F_\varepsilon (S^n)$: 
$$\left \{\matrix 
z_r.z_s = \varepsilon^{-1} z_sz_r \;\,\qquad\quad \quad\quad\;\text{if}\quad\;\; r < s\\
z_r. z_s^* = \varepsilon^{-1}z_s^*z_r\;\, \qquad\quad\quad\quad\;\text{if}\quad\; \;r \ne s\\
 z_r. z_r^* - z_r^*. z_r + (\varepsilon^{-2} + 1) \sum\limits_{s > r} z_s. z_s^* =
0,\\ 
 \sum\limits_{s=0}^n z_s. z_s^* =
0.\qquad\qquad\qquad\qquad\qquad\quad\endmatrix\right.\tag * $$ 

\n Hence, $\Cal F_{\varepsilon} (S^n)$ has (*) as its defining relations. The
construction of irreducible * - representations of $\Cal F_\varepsilon (S^n)$,
is given by 
\vb

\n {\bf Theorem 2.6}. ([3], 13.2.9). {\it Every irreducible *-representation of
$\Cal F_\varepsilon (S^n)$ is equivalent exactly to one of the following: 

i) the one-dimensional representation $\rho_{0, t}, t \in S^1$, given by
$\rho_{0, t} (z_0^*) = t^{-1}, \rho_{0, t} (z_r^*) = 0$ if $r>0$. 

ii) the representation $\rho_{r, t}$, $1 \le r \le n$, $t \in S^1$, on the
Hilbert space tensor product $\ell^2 (\N)^{\otimes r}$, given by 
$$\rho_{r,t} (z_s^*) (e_{k_1} \otimes \cdots \otimes e_{k_r}) =$$
\vskip-0.5cm

$$\cases & \varepsilon^{-(k_1 + \cdots k_s +s)}
 (1-\varepsilon^{-2(k_{s+1}+1)})^  {1/2}  e_{k_1} \otimes \cdots \otimes
 e_{k_s} \otimes e_{ks+1} +1 \otimes e_{ks+2}\\
&\qquad \qquad\qquad\qquad\qquad\hskip1.5cm \otimes \cdots \otimes e_{ks} \;\,
\text{if}\; \;\,s<r\\
&t^{-1}. \varepsilon^{-(k_1+\cdots + k_r + r)} e_{k_1} \otimes \cdots
e_{k_r}\qquad\qquad \hskip0.5cm \;\;
\text{if}\;\;\, r = s\\
&0 \qquad\qquad\qquad\qquad\quad \hskip3.5cm \text{if}\;\;\, s>r \endcases $$

The representation $\rho_{0, t}$ is equivalent to the restriction of the
representation $\Cal T_t$ of $\Cal F_\varepsilon (SL_{n+1} (\C))$ (cf. 2.3);
and for $r>0, \rho_{r, t}$ is equivalent to the restriction of $\pi_{s_1}
\otimes \cdots \otimes \pi_{s_r} \otimes \Cal T_t$.} \hfill \v

From Theorem 2.6, we have 

$$\bigcap\limits_{(w, t) \in W \times\bold T} \text{ker}\; \rho_{w, t} = \{
e\},$$ 
i.e. the representation $\bigoplus\limits_{w \in W} \int_T^\oplus
\rho_{w, t} dt$ is faithful and 
$$\dim \rho_{w, t} = \cases 1 & \text{if}\;\, w = e\\
0 & \text  {if}\;\, w \ne e.\endcases$$ 
We recall now the definition of compact quantum of
spheres $C^*$-algebra.
\vb

\n {\bf Definition 2.7}. {\it The $C^*$-algebraic compact quantum sphere
$C_\varepsilon^* (S^n)$ is the $C^*$-completion of the *-algebra $\Cal
F_\varepsilon (S^n)$ with respect to the $C^*$-norm
$$\Vert f\Vert = \sup\limits_\rho \Vert \rho (f)\Vert\;\;\; (f \in \Cal
F_\varepsilon (S^n)),$$ 
where $\rho$ runs through the *-representations of $\Cal F_\varepsilon (S^n)$
(cf., Theorem 2.6) and the norm on the right-hand side is the operator norm. }

It suffcies to show that $\Vert f\Vert$ is finite for all $f \in \Cal
F_\varepsilon (S^n)$, for it is clear that $\Vert \,.\,\Vert$ is a $C^*$-norm,
i.e. $\Vert f.f^*\Vert = \Vert f\Vert^2$. We now prove the following result
about the structure of compact quantum $C^*$-algebra of sphere $S^n$.
\vb

\n {\bf Theorem 2.8}. {\it With notation as above, we have }
$$C_\varepsilon^* \cong \C(S^1) \oplus \bigoplus\limits_{e\ne w\in W}
\int_{S^1}^\oplus \Cal K (H_{w, t}) dt,$$
{\it where $\text \C(S^1)$ is the algebra of complex valued continuous functions on $S^1$
and $\Cal K (H)$ the ideal of compact operators in a separable Hilbert space
$H$. }
\vb

\n {\it Proof}: Let $w = s_{i_1}. s_{i_2}\ldots s_{i_k}$ be a reduced
decomposition of the element $w \in W$ into a product of reflections. Then by
Proposition 2.6, for $r>0$, the representation $\rho_{w, t}$ is equivalent to the
restriction of $\pi_{s_{i_1}} \otimes\pi_{s_{i_2}} \otimes \cdots \otimes
\pi_{s_{i_k}} \otimes \Cal T_t$, where $\pi_{s_i}$ is the composition of the
homomorphism of $\Cal F_\varepsilon (G)$ onto $\Cal F_\varepsilon (SL_2(\C))$
and the representation $\pi_{-1}$ of $\Cal F_\varepsilon (SL_2(\C))$ in the
Hilbert space $\ell^2(\N)^{\otimes r}$; and the family of one-dimensional
representations $\Cal T_t$, given by 
$$\Cal T_t (a) = t,\;\, \Cal T_t (b) = \Cal T_t (c) = 0,\;\, \Cal T_t (d) =
t^{-1},$$ 
where $t \in S^1$ and $a, b, c, d$ are given by: Algebra $\Cal F_{\varepsilon}
(SL_2(C))$ is generated by the matrix elements of type $A = \pmatrix a & b\\
c & d\endpmatrix$. 

\n Hence, by construction, the representation $\rho_{w, t} = \pi_{s_{i_1}} \otimes
\pi_{s_{i_2}} \otimes \cdots \otimes \pi_{s_{i_k}} \otimes \Cal T_t$. Thus, we
have
\vskip-0.3cm
 
$$\pi_{si}:\, C_\varepsilon^* (S^n) \longrightarrow C_\varepsilon^* (SL_2(\C))
\overset {\pi_{-1}} \to \longrightarrow \Cal L (\ell^2 (\N)^{\otimes r}).$$ 
Now, $\pi_{s_i}$ is CCR (see, [11]) and so, we have $\pi_{s_i}
\big(C_\varepsilon^*(S^n)\big) \cong \Cal K (H_{w, t})$. Moreover $\Cal
T_t(C_{\varepsilon}^* (S^n)) \cong \C$. 

\n Hence,
\vskip-0.8cm

$$\align \rho_{w, t} \big(C_\varepsilon^* (S^n)\big) & = \big(\pi_{s_{i_1}}
\otimes \cdots \otimes \pi_{s_{i_k}} \otimes\Cal T_t\big) \big(C_\varepsilon^* (S^n)\big)\\
& = \pi_{s_{i_1}} \big(C_\varepsilon^* (S^n)\big) \otimes \cdots \otimes
\pi_{s_{i_k}} \big(C_\varepsilon^* (S^n)\big) \otimes \Cal T_t
\big(C_\varepsilon^* (S^n)\big)\\
& \cong \Cal K (H_{s_{i_1}}) \otimes \cdots \otimes \Cal K (H_{s_{i_k}})
\otimes \C\\
& \cong \Cal K (H_{w, t}),\endalign$$ 
\vskip-0.4cm
where \qquad \;\;\, $H_{w, t} = H_{s_{i_1}} \otimes \cdots \otimes H_{s_{i_k}} \otimes \C$. 

\n Thus, \qquad \;\;\, $\rho_{w, t} \big(C_\varepsilon^* (S^n)\big) \cong \Cal K(H_{w, t})$. 

\n Hence, \qquad\quad\; $\bigoplus\limits_{w \in W} \int_{S^1}^\oplus
\rho_{w, t} \big(C_\varepsilon^* (S^n)\big) \cong \bigoplus\limits_{w \in W}
\int_{S^1}^\oplus \Cal K (H_{w, t}) dt.$ 

Now, recall a result of S. Sakai's from [11]: Let $A$ be a commutative
$C^*$-algebra and $B$ be a $C^*$-algebra. Then, $C_0(\Omega, B) \cong A \otimes
B$, where $\Omega$ is the spectrum space of $A$. 

Applying this result, for
$ B  = \Cal K (H_{w, t}) \cong \Cal K$ and
$A  = \C (W \times S^1)$
be a commutative $C^*$-algebra. 

\n Thus, we have
$$C_\varepsilon^* (S^n) \cong \C (S^1) \oplus \bigoplus\limits_{e\ne w \in
W} \int_{S^1}^\oplus \Cal K (H_{w, t}) dt.$$ 
\vskip-0.5cm \hfill\v

Now, we first compute the $K_* \big(C_\varepsilon^* (S^n)\big)$ and the $HE_*
(C_\varepsilon^* (S^n))$ of $C^*$-algebra of quantum sphere $S^n$. 
\vb

\n {\bf Propsition 2.9}. 
\vskip-0.5cm 

$$HE_* \big(C_\varepsilon^*(S^n)\big) \cong H_{DR}^* (W\times S^1).$$
\vb

\n {\it Proof.} We have 
\vskip-0.8cm

$$\align & HE_*\big(C^*_{\varepsilon} (S^n)\big) = 
 HE_*\big(\C(S^1)\oplus \bigoplus\limits_{e\ne w \in W}
\int_{S^1}^\oplus \Cal K (H_{w, t}) dt)\\
& \cong HE_*\big(\C(S^1)\big) \oplus HE_*\Big(\bigoplus\limits_{e\ne w \in W}
\int_{S^1}^\oplus \Cal K (H_{w, t}) dt\Big)\\
& \cong HE_*\big(\C(W\times S^1)\otimes \Cal K\big)\quad\quad (\text{by Proposition 1.1}\; \S1)\\
& \cong HE_*\big(\C(W\times S^1)\big)\endalign$$

Since $C(W\times S^1)$ is a commutative $C^*$-algebra, by Proposition
1.5, \S1, we have
$$ HE_* \big(C_\varepsilon^* (S^n)\big) \cong HE_* \big(\C(W\times S^1)\big) 
\cong H_{DR}^* (W\times S^1)$$
\vb

\n {\bf Proposition 2.10}. 
$$K_* \big(C^*_\varepsilon (S^n)\big) \cong K^* (W\times S^1).$$ 
\vb

\n {\it Proof.} We have
$$\align K_*\big(C^*_\varepsilon (S^n)\big) &= 
K_*\big(\C(S^1)\oplus \bigoplus\limits_{e\ne w \in W}
\int_{S^1}^\oplus \Cal K (H_{w, t}) dt\big)\\
&\cong K_*\big(\C(S^1)\big) \oplus K_*\Big(\bigoplus\limits_{e\ne w \in W}
\int_{S^1}^\oplus \Cal K (H_{w, t}) dt\Big)\\
& \cong K_* (\C(W\times S^1)\otimes \Cal K\big)\\
& \cong K_*\big(\C(W\times S^1)\big)\quad\quad\, (\text{by proposition 1.1}\; \S1).
\endalign$$

In virtue of Proposition 1.5, \S1, we have
$$ K_*\big(\C(W \times S^1)\big) \cong 
K^* (W \times S^1)$$
\vskip-0.5cm \hfill \v
\vb

\n {\bf Theorem 2.11}. {\it With notation above, the Chern character of
$C^*$-algebra of quantum sphere $C_\varepsilon^* (S^n)$
$$ch_{C^*}\,:\, K_* \big(C_\varepsilon^* (S^n)\big) \longrightarrow HE_*
\big(C_\varepsilon^* (S^n)\big)$$ 
is an isomorphism. }
\vb

\n {\it Proof.} By Proposition 2.9 and 2.10, we have: 
$$\align HE_* \big(C_\varepsilon^* (S^n)\big) & \cong HE_* \big(\C(W\times S^1)\big)
\cong H_{DR}^* (W\times S^1),\\
K_* \big(C_\varepsilon^* (S^n)\big) &\cong K_* \big(\C (W\times S^1)\big) \cong K^* (W\times S^1).\endalign$$

Now, consider the commutative diagram
$$\matrix K_* \big(C_\varepsilon^* (S^n)\big) & \overset{ch_{C_\varepsilon^*}} \to
{\text{-----}\hskip-0.2cm 
\longrightarrow } &
HE_* \big(C^*_\varepsilon(S^n)\big)\\
\bigg\downarrow \cong && \bigg\downarrow \cong \\
K_* \big(\C(W\times S^1)\big) & \overset {ch_{CQ}}
\to{\text{-----}\hskip-0.2cm \longrightarrow } & 
HE_* \big(\C(W\times S^1)\big)\\ 
\bigg\downarrow \cong && \bigg\downarrow \cong \\
K^* \big(W\times S^1\big) & \overset{ch}\to {\text{-----}\hskip-0.2cm
\longrightarrow} & H_{DR}^* 
\big(W\times S^1\big)\endmatrix $$  

Moreover, following Watanabe [15], the $ch\,:\,K^* (W\times S^1) \otimes \C
\longrightarrow H_{DR}^* (W\times S^1)$ is an isomorphism. 

Thus, \quad $ch_{C_\varepsilon^*}\,:\,K^* \big(C_\varepsilon^* (S^n)\big)
\longrightarrow HE_* \big(C_\varepsilon^*(S^n)\big)$ in an isomorphism. 
\vb

\n {\bf Acknowledgment.} {\it The author would like to thank Professor Do Ngoc Diep
for his guidance and encouragement during this work. }
\vskip0.7cm

{\centerline\bf REFERENCES}
\vl

\item {[1]. } J. Cuntz, {\it Entice cyclic cohomology of Banach algebra and character
of $\theta$-summable Fredhom modules}, $K$-Theory, {\bf 1}(1998), 519-548. 
\vn

\item {[2]. } J. Cuntz and D. Quillen, {\it The $X$ complex of the unuversal
extensions,} 
Preprint Math. Inst. Uni. Heidelbeg, 1993. 
\vn

\item {[3]. } V. Chari and A. Pressley, {\it A guide to Quantum groups}, Cambridge
Uni. Press, 1995. 
\vn

\item {[4]. } D.\, N. Diep, A.\, O. Kuku and N.\, Q. Tho, {\it Non-commutative Chern
character of compact Lie group $C^*$-algebras}, $K$- Theory, {\bf 17(2)}, (July 1999),
195-208. 
\vn

\item {[5]. } D.\, N. Diep, A.\, O. Kuku and N.\, Q. Tho, {\it Non-commutative Chern
character of compact quantum groups }, ICTP (Preprint), IC/98/91, 1998, 16pp. 
\vn

\item {[6]. }D.\, N. Diep and N.\, V. Thu, {\it Homotopy invariance of entire curent
cyclic homology}, Vietnam J. of Math., {\bf 25(3)} (1997), 211-228. 
\vn

\item {[7]. } D.\, N. Diep and N.\, V. Thu, {\it Entire homology of non-commutative
de.Rham curents}, ICTP, IC/96/214, 1996, 23pp; to apprear in Publication of
CFCA, Hanoi-City Vietnam National Iniversity, 1997. 
\vn

\item {[8]. }M. Khalkhali, {\it On the intive cyclic cohomology of Banach algebras:
I. Morita invariance}, Mathematisches Inst. Uni. Heidelberg, {\bf 54} (1992), pp24. 
\vn

\item {[9]. } M. Khalkhali, {\it On the entite cyclic cohomology of Banach algebras:
II. Homotopy invariance}, Mathematisches Inst. Uni. Heidelberg, {\bf 55} (1992), pp18.
\vn

\item {[10]. }J. Packer, {\it Transformation group $C^*$-algebra: A selective
survey}, contemporary Mathematics, Volume {\bf 167}, 1994. 
\vn

\item {[11]. } S. Sakai, $C^*${\it -algebras and $W^*$-algebras}, Spriner-Verlag
Berlin. Heidelberg New York, 1971. 
\vn

\item {[12]. }N.\, V. Thu, {\it Morita invariance of entire current cyclic homology},
Vietnam J. Math. {\bf 26}: 2(1998), 177-179. 
\vn

\item {[13]. } N.\, N. Thu, {\it Exactness of entire current cyclic homology}, Acta
Math, Vietnamica (to appear). 
\vn

\item {[14]. } T. Watanabe, {\it Chern character on compact Lie groups of low rank},
Osaka J. Math. {\bf 2} (1985), 463-488. 
\vn

\item {[15]. } T. Watanabe, {\it On the Chern character of symmetric space related
to} $SU (n)$, J. Math. Kyoto Uni. {\bf 34} (1994), 149-169. 
\vb


\enddocument